\documentclass{amsart}

\usepackage{graphicx}
\usepackage{tikz}
\usepackage{tikz-cd}
\usepackage{verbatim}

\newtheorem{theorem}{Theorem}[section]

\theoremstyle{definition}
\newtheorem{definition}[theorem]{Definition}
\newtheorem{example}[theorem]{Example}

\theoremstyle{remark}

\theoremstyle{proposition}
\newtheorem{proposition}[theorem]{Proposition}

\theoremstyle{corollary}

\numberwithin{equation}{section}

\DeclareMathOperator{\Hom}{Hom} 
 \DeclareMathOperator{\ann}{ann}

\begin{document}

\title{Absolutely Self Pure Modules}


\thanks{}

\author{Mohanad Farhan Hamid}

\address{Department of Mathematics, College of Education, Misan University, Amarah, Iraq}
\curraddr{} \email{mohanad@uomisan.edu.iq}

\subjclass[2000]{16D50}
\keywords{(Self) pure submodule, absolutely (self) pure module, quasi injective module}

\begin{abstract} An $R$-module $M$ is called absolutely self pure if for any finitely generated left ideal of $R$ whose kernel is in the filter generated by the set of all left ideals $L$ of $R$ with $L \supseteq \ann (m)$ for some $m \in M$, any map from $L$ to $M$ is a restriction of a map $R \rightarrow M$. Certain properties of quasi injective and absolutely pure modules are extended to absolute self purity. Regular and left noetherian rings are characterized using this new concept.
\end{abstract}

\maketitle

\section{}
Let $R$ be an associative ring with identity. By a module and an $R$-module we mean a left unital $R$-module. A submodule $A$ of an $R$-module $B$ is called a \emph{pure submodule} \cite{Sten} if for any commutative diagram with $K$ a finitely generated submodule of a free $R$-module $F$
$$\begin{tikzcd}
K \arrow[hook]{r}\arrow{d}{}
&F\arrow{d}{} \arrow[dashed]{ld}\\
A\arrow[hook]{r}&B
\end{tikzcd}$$
there is a map $F \rightarrow A$ making the upper triangle
commute.
A module $A$ is called \emph{absolutely pure} \cite{Mdx} if it is pure in every module containing it as a submodule or equivalently if $A$ is pure in some injective module or $A$ is pure in its injective envelope \cite{M}. Our goal is to study modules that are pure in their quasi injective envelopes. Of course, if the module is pure in every quasi injective module then it must be absolutely pure. So we give a weaker form of purity, \emph{self purity}, and study modules that are self pure in every module containing them as submodules. This turns out to be equivalent to saying that the module is self pure in some quasi injective module or self pure in its quasi injective envelope.

The quasi injective envelope of an $R$-module $M$ is denoted $Q(M)$. By $\Omega(M)$ we mean the set of all left ideals $L$ of $R$ such that $L \supseteq \ann (m)$ for some $m \in M$. The filter generated by $\Omega(M)$ in the lattice of left ideals of $R$ is denoted $\overline{\Omega}(M)$. Recall that L. Fuchs \cite{Fuchs} has proved that a module $M$ is quasi injective if and only if for any homomorphism $f$ from a left ideal $L$ of $R$ whose kernel is in $\overline{\Omega}(M)$ there is an extension to a homomorphism $R \rightarrow M$ of $f$. This gives a generalization of Baer condition for injective modules.

On the other hand, absolutely pure modules $A$ are characterized by the property that for any finitely generated submodule $K$  of a free module $F$, any homomorphism $K \rightarrow A$ can be extended to a homomorphism $F \rightarrow A$ \cite{M}.

Generalizing the above facts, we say that a module $A$ is absolutely self pure if any map from a finitely generated left ideal of $R$ into $A$ whose kernel is in $\overline{\Omega}(A)$ can be extended to a map $R \rightarrow A$. Regular and left noetherian rings are characterized using properties of absolutely self pure modules.
\section{}
\begin{definition}
A submodule $A$ of an $R$-module $B$ is called a \emph{self pure}
submodule of $B$ (denoted $A \leq^{sp} B$) if the following condition holds:
For any finitely
generated left ideal $L$ of $R$ and any map $f : L \rightarrow
A$ with $\ker f \in \overline{\Omega}(A)$, if there is a
map $R \rightarrow B$ making the following diagram
commutative:
$$\begin{tikzcd}
L \arrow[hook]{r}\arrow{d}{f}
&R\arrow{d} \arrow[dashed]{ld}\\
A\arrow[hook]{r}&B
\end{tikzcd}$$
then there exists a map $ R \rightarrow A$ making the upper triangle
commutative.
\end{definition}

Any pure submodule of a module is self pure, but not conversely,
for one may take any quasi injective module that is not absolutely
pure, which must clearly be self pure in its injective envelope,
but of course not pure.

More generally, one can study purity with respect to another module $M$. Precisely, we say that a submodule $A$ of an $R$-module $B$ is called $M$-\emph{pure}
submodule of $B$ if the following condition holds:
For any finitely
generated left ideal $L$ of $R$ and any map $f : L \rightarrow
A$ with $\ker f \in \overline{\Omega}(M)$, if there is a
map $R \rightarrow B$ making the above diagram
commutative then there exists a map $ R \rightarrow A$ making the upper triangle
commutative. Therefore, $A$ is $A$-pure in $B$ exactly when $A$ is self pure in $B$. However we will restrict our attention to the concept of self purity.
\begin{proposition} \label{sptrans}
If $A$, $B$ and $C$ are $R$-modules such that $A \leq^{sp} B$ and
$B \leq^{sp} C$ then $A \leq^{sp} C$.
\begin{proof}
Consider the following diagram
$$\begin{tikzcd}
L \arrow[hook]{r}\arrow{d}{f}
&R\arrow{d} \\
A\arrow[hook]{r}&C
\end{tikzcd}$$
where $L$ is a finitely generated left ideal of $R$ with $\ker f \in
\overline{\Omega}(A)$ and $R \rightarrow C$ is a map making the
diagram commutative. Then obviously, we can consider $f$ as a map
$L \rightarrow B$ for which we have a commutative diagram:
$$\begin{tikzcd}
L \arrow[hook]{r}\arrow{d}{f}
&R\arrow{d} \\
B\arrow[hook]{r}&C
\end{tikzcd}$$
Since $B \leq^{sp} C$, there is an extension $g : R \rightarrow B$ of $f$, so that the following square commutes:
$$\begin{tikzcd}
L \arrow[hook]{r}\arrow{d}{f}
&R\arrow{d}{g} \\
A\arrow[hook]{r}&B
\end{tikzcd}$$
and since $A \leq^{sp} B$, there is an $h : R \rightarrow A$ extending $f$.
\end{proof}
\end{proposition}

\begin{proposition} \label{spantitrans}
If $A \subseteq B \subseteq C$ are $R$-modules such that $A \leq^{sp} C$ then $A \leq^{sp} B$.
\begin{proof}
Consider the following commutative diagram with the obvious maps:
$$\begin{tikzcd}
L \arrow[hook]{r}\arrow{d}{f}
&R\arrow{d} &\\
A\arrow[hook]{r}&B \arrow[hook]{r}&C
\end{tikzcd}$$
Since $A \leq^{sp} C$, there is a $g : R \rightarrow A$ extending $f$.
\end{proof}
\end{proposition}

\begin{definition}
A module is called \emph{absolutely self pure} if it is self pure in every module containing it as a submodule.
\end{definition}

It is clear that absolutely pure modules and quasi injective modules are examples of absolutely self pure modules. If a module $A$ is self pure in some quasi injective module then by Proposition \ref{spantitrans}, it must be self pure in its quasi injective envelope $Q(A)$. So if $A$ is contained in some other module $B$, then as $Q(A) \leq^{sp} Q(B)$ we must have by Propositions \ref{sptrans} and \ref{spantitrans} that $A \leq^{sp} B$. Therefore, a module is absolutely self pure if and only if it is self pure in some quasi injective module if and only if it is self pure in its quasi injective envelope.
\begin{theorem} \label{absspurechar}
A module $A$ is absolutely self pure if and only if for each finitely generated left ideal $L$ of $R$ and each map $f : L \rightarrow A$ with $\ker f \in \overline{\Omega}(A)$ there is an extension map $ R \rightarrow A$ of $f$.
\begin{proof}
Consider the following commutative diagram:
$$\begin{tikzcd}
L \arrow[hook]{r}\arrow{d}{f}
&R\arrow{d}{g} &\\
A\arrow[hook]{r}&Q(A)
\end{tikzcd}$$
where $Q(A)$ denotes the quasi injective envelope of $A$. Existence of $g$ is guaranteed be quasi injectivity of $Q(A)$. Now $A$ is self pure in $Q(A)$ if and only if $f$ can be extended to a map $R \rightarrow A$.
\end{proof}
\end{theorem}

By Proposition \ref{sptrans}, a self pure submodule of an absolutely self pure module is again absolutely self pure. In particular, direct summands of absolutely self pure modules are absolutely self pure.

\begin{theorem}
A module $A$ is absolutely self pure if and only if any direct sum of copies of $A$ is absolutely self pure.
\begin{proof}
Suppose that $A$ is absolutely self pure and let $f$ be a map from a finitely generated left ideal $L$ of $R$ into $A^{(I)}$ for some index set $I$ such that $\ker f \in \overline{\Omega}(A^{(I)})$, i.e. $\ker f \supseteq \bigcap_{j \in J} \ann ((a_i)_j)$ for finite $J$. But the non-zero $a_i$'s in each $(a_i)_j$ are also finite. Therefore $\ker f$ contains the intersection of annihilators of a finite set of individual $a_i$'s and hence it belongs to $\overline{\Omega}(A)$.
Let $\{l_1, \cdots, l_n\}$ be a generating set for $L$.  Each of the images $f(l_i)$ is of finite support, so $f$ can be considered as a map $L \rightarrow A^{(K)}$ for some finite subset $K$ of $I$ and hence it is the finite direct sum of coordinate maps $f_k$ into each factor $A$. Clearly $\ker f_k \supseteq \ker f$ and so $\ker f_k \in \overline{\Omega}(A)$ for each coordinate map $f_k$. By absolute self purity of $A$, each $f_k$ is extendable to a $g_k : R \rightarrow A$. The map whose coordinate maps are the $g_k$'s is the desired extension of $f$. The other implication is clear.
\end{proof}
\end{theorem}

Left noetherian rings are precisely rings over which every absolutely pure module is (quasi) injective \cite[Theorem 3]{M} and \cite[Theorem 8]{(Flat)modules}. It is clear that over such rings, the concepts of absolutely
self pure and that of quasi injective modules are
equivalent and by \cite[Theorem 8]{(Flat)modules} the converse is also true. So we have the following characterization of left
noetherian rings:
\begin{theorem} \label{noethchar}
A ring $R$ is left noetherian if and only if any absolutely self
pure $R$-module is quasi injective.
\qed
\end{theorem}

Recall that a ring $R$ is called regular if every principal left ideal of $R$ is a direct summand. Over a regular ring, all modules are absolutely pure and hence all modules are absolutely self pure. The converse is also true:

\begin{theorem} \label{regchar}
A ring $R$ is regular if and only if every left $R$-module is absolutely self pure.
\begin{proof}
($ \Rightarrow $) Clear. ($ \Leftarrow $) Given a principal left ideal $L$ of $R$, extend the identity map of $L$ to the map $L \rightarrow L \oplus R$ defined by $a \mapsto (a,0)$. This map clearly has kernel containing $\ann (0,1)$. By assumption $L \oplus R$ is absolutely self pure, hence there is a $g : R \rightarrow L \oplus R$ extending $f$. Follow $g$ by the projection $L \oplus R \rightarrow L$ to get an extension of the identity map $L \rightarrow L$, which shows that $L$ is a direct summand of $R$.
\end{proof}
\end{theorem}

Combining Theorems \ref{noethchar} and \ref{regchar} we get the well-known characterization that $R$ is semisimple (= regular and left noetherian) if and only if every $R$-module is quasi injective.
 
By Theorem \ref{noethchar}, over the ring of integers $\mathbb Z$ the absolutely self pure modules are exactly the quasi injective ones. So any quasi injective abelian group that is not injective serves as an example of an absolutely self pure module which is not absolutely pure. Also by Theorems \ref{regchar} and \ref{noethchar}, if a ring $R$ is regular but not left noetherian then there must exist an absolutely self pure $R$-module that is not quasi injective. Another example is the following.

\begin{example}
We know that if $R$ is a regular but not a left noetherian ring then there must exist an absolutely pure $R$-module $A$ that is not quasi injective. Over the ring $\mathbb Z$ the module $\mathbb Z_2$ is quasi injective but not absolutely pure. Now let $R$ and $A$ be as above and consider the module 
$M=\left( {\begin{array}{c}A  \\
    \mathbb Z_2 \\
  \end{array} } \right)$ over the ring 
   $S=\left( {\begin{array}{cc}
   R & 0 \\
   0 & \mathbb Z \\
  \end{array} } \right)$
The module $M$ is not quasi injective (nor absolutely pure), for otherwise the direct summand $\left( {\begin{array}{c}A  \\
    0 \\
  \end{array} } \right)$ would be quasi injective (the direct summand $\left( {\begin{array}{c}0  \\
    \mathbb Z_2 \\
  \end{array} } \right)$ would be absolutely pure). Now we proceed to show that $M$ is absolutely self pure. Any finitely generated left ideal of $S$ is of the form $K=\left( {\begin{array}{cc}
   I & 0 \\
   0 & J \\
  \end{array} } \right)$, where $I$ (respectively $J$) is a finitely generated left ideal of $R$ (respectively $\mathbb Z$). Any map $K \rightarrow M$ is of the form $\left({\begin{array}{cc}
   f &0\\
   0&g \\
  \end{array} } \right)$ where $f \in \Hom_R(I,A)$ and $g \in \Hom_{\mathbb Z}(J, \mathbb Z_2)$. To extend $\left({\begin{array}{cc}
   f &0\\
   0&g \\
  \end{array} } \right)$ to $S$ is to extend both $f$ to $R$ and $g$ to $\mathbb Z$. The former is always possible because $I$ is a finitely generated left ideal of a regular ring, hence it is a direct summand. So we focus on extending $g$ to a map $\mathbb Z \rightarrow \mathbb Z_2$. But if $\ker (K \rightarrow M)$ contains the intersection of $\ann \left( {\begin{array}{c}a_i  \\
    z_i \\
  \end{array} } \right)$ for some $\left( {\begin{array}{c}a_i  \\
    z_i \\
  \end{array} } \right) \in M$, $i=1, \cdots, n$, which must be equal to $\left( {\begin{array}{cc} \cap \ann (a_i) & 0  \\
    0&\cap \ann (z_i) \\
  \end{array} } \right)$ then we need only to focus on $\cap \ann (z_i)$. This latter is equal to $2\mathbb Z$ or $\mathbb Z$ according to whether one of the $z_i$'s is $1$ or not. In any case we must have that $2 \mathbb Z \subseteq \ker (g:J \rightarrow \mathbb Z_2)$. Hence either $J = \mathbb Z$ and we are done or $J = 2\mathbb Z$ and therefore $g$ is the zero map.
\end{example}
A note about modules that are pure in their quasi injective envelopes is in order. Let us call them \emph{absolutely quasi pure} modules. It is obvious that this concept lies between quasi injectivity and absolute self purity. It is not clear whether every absolutely self pure module is absolutely quasi pure. Of course over regular rings the two concepts are equivalent, and over left noetherian rings they are also equivalent to quasi injectivity. If a module $A$ is absolutely quasi pure then so is any finite direct sum of copies of $A$.
\bibliographystyle{amsplain}

\end{document}